%% file: main.tex
\documentclass[11pt,letterpaper]{article}
\pdfoutput=1 

\input{preamble.tex}
\input{def.tex}

\usepackage{lipsum}

\newcommand*\samethanks[1][\value{footnote}]{\footnotemark[#1]}

\title{Primal heuristics for Dantzig-Wolfe decomposition for unit~commitment}
\author{Nagisa Sugishita\thanks{School of Mathematics, University of Edinburgh, James Clerk Maxwell Building, Edinburgh (UK), EH9 3FD, Email: {\tt n.sugishita@sms.ed.ac.uk}, {\tt a.grothey@ed.ac.uk}, {\tt k.mckinnon@ed.ac.uk}}, Andreas Grothey\samethanks, Ken McKinnon\samethanks}

\begin{document}

\maketitle

\begin{abstract}
\input{ph_abstract}
\end{abstract}

\glsresetall

\section{Introduction}
\label{ph_sec_introduction}

\input{ph_sec_introduction}

\section{Dantzig-Wolfe decomposition}
\label{ph_sec_dantzig_wolfe_decomposition}

\input{ph_sec_dantzig_wolfe_decomposition}

\section{Primal heuristics}
\label{ph_sec_primal_heuristics}

\input{ph_sec_primal_heuristics}

\section{Numerical experiments}
\label{ph_sec_numerical_experiments}

\input{ph_sec_numerical_experiments}

\section{Conclusion}
\label{ph_sec_conclusion}

\input{ph_sec_conclusion}


\printbibliography

\input{ph_appendix}


\end{document}

%% file: preamble.tex
\usepackage{amsmath, amssymb, amsthm, amsfonts}
\usepackage[font=small,labelfont=bf]{caption}
\usepackage{graphicx}
\usepackage[outdir=./]{epstopdf}
\usepackage{xcolor}
\usepackage{booktabs}
\usepackage{algorithm}
\usepackage{algpseudocode}

\usepackage[xindy,acronyms]{glossaries}
\glsdisablehyper  

\usepackage{csvsimple}  
\usepackage{tikz}
\usetikzlibrary{arrows,calc,matrix,patterns,positioning}

\usepackage{hyperref}
\hypersetup{colorlinks=true,
            urlcolor=blue,
            linkcolor=black,
            citecolor=blue,
            bookmarksdepth=paragraph,
            pdftitle={Primal heuristics for Dantzig-Wolfe decomposition for unit commitment},
            pdfauthor={Nagisa Sugishita, Andreas Grothey, Ken McKinnon},}

\usepackage{pgfplots, pgfplotstable}
\pgfplotsset{compat=1.17}
\usepgfplotslibrary{statistics}
\usepgfplotslibrary{groupplots}
\usepgfplotslibrary{dateplot}
\usepgfplotslibrary{colorbrewer}

\usepackage[
  style=numeric,
  backend=biber,
  maxcitenames=2,
]{biblatex}

\theoremstyle{plain}

\theoremstyle{plain}

\theoremstyle{plain}

\theoremstyle{plain}

\theoremstyle{definition}

\theoremstyle{definition}

\theoremstyle{definition}

\theoremstyle{definition}

\theoremstyle{remark}

\numberwithin{equation}{section}

\newacronym{DW}{DW}{Dantzig-Wolfe}
\newacronym{LP}{LP}{linear programming}
\newacronym{LPR}{LPR}{linear programming relaxation}
\newacronym{MILP}{MILP}{mixed-integer linear programming}
\newacronym{MP}{MP}{master problem}
\newacronym{RMP}{RMP}{restricted master problem}
\newacronym{UC}{UC}{unit commitment}

\newcommand{\citet}[1]{\textcite{#1}}

%% file: def.tex
\definecolor{mysymbolcolor}{rgb}{0,0,0}

\providecommand\idxcomponent{}
\renewcommand{\idxcomponent}{{\color{mysymbolcolor} g}}

\providecommand\numcomponents{}
\renewcommand{\numcomponents}{{\color{mysymbolcolor} G}}

\providecommand\oneofidxcomponents{}
\renewcommand{\oneofidxcomponents}{%
{\color{mysymbolcolor} \idxcomponent = 1, 2, \ldots, \numcomponents}%
}

\providecommand\sumcomponents{}
\renewcommand{\sumcomponents}[1][\idxcomponent]{%
{\color{mysymbolcolor} \sum_{#1 = 1}^\numcomponents}%
}

\providecommand\dualfunction{}
\renewcommand{\dualfunction}{{\color{mysymbolcolor} q}}

%% file: ph_abstract.tex
The unit commitment problem is a short-term planning problem in the energy industry.
Dantzig-Wolfe decomposition is a popular approach to solve the problem.
This paper focuses on primal heuristics used with Dantzig-Wolfe decomposition.
We propose two primal heuristics: one based on decomposition and one based on machine learning.
The first one uses the fractional solution to the restricted master problem to fix a subset of the integer variables.
In each iteration of the column generation procedure, the primal heuristic obtains the fractional solution, checks whether each binary variable satisfies the integrality constraint and fix those which do.
The remaining variables are then optimised quickly by a solver to find a feasible, near-optimal solution to the original instance.
The other primal heuristic based on machine learning is of interest when the problems are to be solved repeatedly with different demand data but with the same problem structure.
The primal heuristic uses a pre-trained neural network to fix a subset of the integer variables.
In the training phase, a neural network is trained to predict for any demand data and for each binary variable how likely it is that the variable takes each of two possible values.
After the training, given an instance to be solved, the prediction of the model is used with a rounding threshold to fix some binary variables.
Our numerical experiments compare our methods with solving the undecomposed problem and also with other primal heuristics from the literature.
The experiments reveal that the primal heuristic based on machine learning is superior when the suboptimality tolerance is relatively large, such as 0.5\% or 0.25\%, while the decomposition is the best when the tolerance is small, for example 0.1\%.

%% file: ph_sec_introduction.tex
The \gls{UC} problem is a combinatorial optimisation problem to find the optimal operating schedule of power plants for given demand over a fixed period.
This problem has been actively studied for the past few decades, see \cite{VanAckooijetal2018a} for a recent survey.

One popular solution method for a \gls{UC} problem is to use Dantzig-Wolfe decomposition to decompose the problem by generators.
The reformulated problem is then solved with a column generation procedure, which yields a sequence of (hopefully improving) lower bounds.
This procedure is the dual of a cutting-plane approach to Lagrangian relaxation, as discussed in \cite{Briantetal2008}.
For problems with integer variables, Dantzig-Wolfe reformulation is a relaxation of the original problem, however \citet{Bertsekasetal1983} and \citet{Bard1988} reported that the optimality gap introduced by the reformulation is typically small, if the problem size is large.

Since Dantzig-Wolfe decomposition only provides lower bounds, it is necessary to be combined with a primal heuristic to obtain a feasible solution.
Typically primal heuristics based on decomposition are used \cite{MerlinAndSnadrin1983,ZhuangAndGaliana1988,Guanetal1992,Takritietal2000}.
Such primal heuristics collect information through the execution of the column generation procedure and use it as a guide to create primal feasible solutions.
Empirically it has been observed that these methods are capable of finding a primal feasible solution with small suboptimality.
However, primal heuristics based on decomposition are not the only options in Dantzig-Wolfe decomposition.
A decomposition-free primal heuristic may be advantageous since primal heuristics based on decomposition typically need a few iterations of the column generation procedure to collect enough data to create solutions of high quality.
This requires some computational time, which may be avoided if a good feasible solution is available earlier.

In recent years, there has been a growing interest in the use of machine learning in primal heuristics.
\citet{Belloetal2016}, \citet{Khaliletal2017}, \citet{Nazarietal2018} and \citet{Kooketal2019} applied reinforcement learning on various combinatorial problems, including the travelling salesman problems and the vehicle routine problems.
\citet{Nairetal2018} used reinforcement learning to solve two-stage stochastic programming approximately.
In all of the above studies, a machine learning model was designed to create primal feasible solutions directly.
There are also approaches which use both machine learning and an optimisation solver.
\citet{BertsimasAndStellato2021} and \citet{BertsimasAndStellato2020} proposed a method for mixed-integer convex/quadratic optimisation problems, where a machine learning model was used to simplify the optimisation problems.
They first trained a model to predict the tight (active) constraints and the values of the integer variables with supervised learning.
Given a new instance, the model was used to delete inactive constraints and fix all the integer variables and the reduced model was solved quickly.
\citet{Xavieretal2020} studied the use of a support vector machine on a security-constrained \gls{UC} problem where the model was used to reduce the problem size by fixing a subset of binary variables.
They first trained a model for each binary variable and then selected all models whose accuracy in predicting their values were higher than a prescribed value.
Then given a new instance the selected models were used to fix the corresponding variables and the reduced problem was then solved with an \gls{MILP} solver.
\citet{Wang2021} proposed a similar approach to reduce the problem size of a knapsack problem, which was based on a nearest neighbour method.
Typically these approaches take a longer time compared with the methods which are purely based on machine learning models, however they are likely to provide solutions of better quality.
We note that many of the above approaches focus on finding a primal feasible solution.
They do not compute the suboptimality of the output and may not be able to deliver a solution whose suboptimality is smaller than a prescribed tolerance, if any.

In our study, we propose two primal heuristics: one based on decomposition and one based on machine learning.
The primal heuristic based on decomposition uses the fractional solution to the \gls{RMP} to fix a subset of integer variables.
That is, in each iteration of the column generation procedure, given the fractional solution to the \gls{RMP}, the primal heuristic checks whether each binary variable satisfies the integrality constraint and fix those which do.
The remaining variables are then optimised quickly by an \gls{MILP} solver with the aim of finding a feasible, near-optimal solution to the original instance.

The other primal heuristic based on machine learning is of interest when the problems are to be solved repeatedly with different demand data but with the same problem structure.
For example, a \gls{UC} problem is solved repeatedly on a daily basis by electricity generating companies.
In such a case, typically the structure of the problem is the same and only some of the problem data such as demand is modified.
Under this setup, the primal heuristic can use a pre-trained model instead of the \gls{RMP} to fix a subset of the integer variables.
In the training phase, a model is trained to predict for any demand data and for each binary variable how likely it is that the variable takes each of two possible values.
After the training, given an instance to be solved, the prediction of the model is used with a rounding threshold to fix some binary variables and reduce the problem size.
Similar approaches were used by \citet{Xavieretal2020,Wang2021}.
However, in our approach we adjust the amount of problem size reduction adaptively depend on the instance to be solved to provide high quality solution quickly.
Furthermore, we combine the primal heuristic with Dantzig-Wolfe decomposition, which allows us to find a primal feasible solution with proven suboptimality in a short time.

We note that in the literature a large number of decomposition-free primal heuristics, such as nature-inspired primal heuristics, to solve \gls{UC} have been proposed as well.
For recent survey covering such primal heuristics, see \citet{Saravananetal2013}.
However, many of them have not been tested comprehensively yet.
Often the quality of output solutions such as suboptimality on large-scale instances is not reported.
Since the capability of such methods to deliver primal feasible solutions with small suboptimality (such as 0.1\%) is not yet clear, we do not consider them further in this paper.

The rest of the paper is structured as follows.
Section \ref{ph_sec_dantzig_wolfe_decomposition} briefly introduces Dantzig-Wolfe decomposition and the column generation procedure.
Section \ref{ph_sec_primal_heuristics} reviews primal heuristics found in the literature and presents our primal heuristics based on decomposition and based on machine learning.
The proposed approaches are tested on large scale \gls{UC} problems in Section \ref{ph_sec_numerical_experiments}.
Finally, in Section \ref{ph_sec_conclusion} conclusions of this work are presented.

%% file: ph_sec_dantzig_wolfe_decomposition.tex
In this section, we briefly review Dantzig-Wolfe decomposition and the column generation procedure.
For further background, see \cite{VanderbeckAndSavelsbergh2006}.

Consider the following family of mixed-integer programmes parameterised by $\omega$:
\begin{align}
z(\omega) = \min_{x_1, x_2, \ldots, x_\numcomponents} \ & \sumcomponents c_\idxcomponent^T x_\idxcomponent
\label{problem_eq_mip} \\
\text{s.t.} \ &
\sumcomponents A_\idxcomponent x_\idxcomponent = a(\omega), \notag \\
& x_\idxcomponent \in X_\idxcomponent(\omega),
\quad \oneofidxcomponents, \notag
\end{align}
where $x_\idxcomponent$ is a vector of decision variables for each $\oneofidxcomponents$, and
\[
X_\idxcomponent(\omega) = \{x_\idxcomponent = (u_\idxcomponent, z_\idxcomponent) \in \mathbb{R}^{n} \times \{0, 1\}^{m} \mid D x_\idxcomponent \le d(\omega)\},
\quad \oneofidxcomponents.
\]
For simplicity, we assume that $X_\idxcomponent(\omega)$ is non-empty and bounded for every $\omega$ and $\idxcomponent$ and that the problem \eqref{problem_eq_mip} has an feasible solution for every $\omega$.
To reduce clutter in what follows we drop the dependence on $\omega$ unless it causes confusion.

In Dantzig-Wolfe decomposition, we consider a relaxation of \eqref{problem_eq_mip}, referred to as the \gls{MP}, got by replacing $X_\idxcomponent$ with $\mathrm{conv}(X_\idxcomponent)$ for every $\idxcomponent$.
Let $\{x_{\idxcomponent i} \mid i \in I_\idxcomponent\}$ denote the extreme points of $X_\idxcomponent$.
Given the boundedness assumption on $X_\idxcomponent$, it follows that
\begin{equation*}
\mathrm{conv}(X_\idxcomponent) =
\left\{
\sum_{i \in I_\idxcomponent} x_{\idxcomponent i} p_{\idxcomponent i} \ \Big\vert \ \sum_{i \in I_\idxcomponent} p_{\idxcomponent i} = 1, \quad p_{\idxcomponent i} \ge 0 \quad (i \in I_\idxcomponent)
\right\},\quad \oneofidxcomponents.
\end{equation*}
This implies that the \gls{MP} can be written as a \gls{LP} with decision variables $\{ p_{\idxcomponent i} \mid \oneofidxcomponents, \ i \in I_\idxcomponent \}$.
However, finding all the extreme points is time consuming and leads to a formulation that is too large to solve explicitly, so a column generation procedure is used.
The \gls{RMP} is defined by replacing $I_\idxcomponent$ in the \gls{MP} with a subset $\hat{I}_\idxcomponent \subset I_\idxcomponent$ for every $\idxcomponent$.
Thus, the \gls{RMP} is given by
\begin{align}
\min_{p} \ & \sumcomponents \sum_{i \in \hat{I}_\idxcomponent} c_\idxcomponent^T x_{\idxcomponent i} p_{\idxcomponent i}
\label{problem_eq_rmp} \\
\text{s.t.} \ &
\sumcomponents \sum_{i \in \hat{I}_\idxcomponent} A_\idxcomponent x_{\idxcomponent i} p_{\idxcomponent i} = a, \notag \\
& \sum_{i \in \hat{I}_\idxcomponent} p_{\idxcomponent i} = 1, \quad \oneofidxcomponents, \notag \\
& p_{\idxcomponent i} \ge 0, \quad \oneofidxcomponents, i \in \hat{I}_\idxcomponent. \notag
\end{align}
Suppose that the \gls{RMP} is feasible (the $\hat{I}_\idxcomponent$'s must have enough elements so that the first constraint in \eqref{problem_eq_rmp} can be satisfied) and let $y$ and $\sigma_\idxcomponent$ for $\oneofidxcomponents$ be the optimal dual solution to the \gls{RMP} corresponding to the first and second constraints respectively.
To find columns to be added to the \gls{RMP}, the following problems, known as the pricing subproblems, are solved:
\begin{equation}
r_\idxcomponent(y) = \min_{x_\idxcomponent} \{(c_\idxcomponent - y^T A_\idxcomponent) x_\idxcomponent \mid x_\idxcomponent \in X_\idxcomponent\},
\quad \oneofidxcomponents.
\label{eq_pricing_problem}
\end{equation}
If $r_\idxcomponent(y) \ge \sigma_\idxcomponent$ for all $\idxcomponent$, the \gls{RMP} has found the optimal solution to the \gls{MP}.
Otherwise, the solutions to the pricing subproblems are added to the set $\hat{I}_\idxcomponent$ and the above process is repeated.
It follows from \gls{LP} duality that given a dual value $y$
\begin{equation}
\dualfunction(y) = a^T y + \sumcomponents r_\idxcomponent(y)
\label{eq_lower_bound}
\end{equation}
is a lower bound to the \gls{MP} and so \eqref{problem_eq_mip}.

It is known that the naive column generation approach described above suffers from instability.
To mitigate the issue, quadratic regularisation on the dual variable may be added to the \gls{RMP}, as described by \citet{Briantetal2008}.
It has been numerically observed that by initialising the regularisation centre to values close to the optimal one we can warmstart the algorithm and reduce the computational time.
\citet{Schulzeetal2017} solved the \gls{LPR} of the original problem \eqref{problem_eq_mip} and used the optimal dual values to the \gls{LPR} to initialise the column generation procedure.
We follow their approach in this study.
Since the column generation procedure only provides lower bounds, we need to run a primal heuristic to obtain a feasible solution, which is discussed in the next section.

The above discussions on the column generation procedure are put together in Algorithm \ref{alg_plain_dw}.
In each iteration, the pricing subproblems are solved and a lower bound is evaluated.
It is followed by \gls{RMP} update, its solution and a primal heuristic.
A more detailed description of the algorithm is found in \citet{Schulzeetal2017}.

\begin{algorithm}
\caption{Column generation procedure}
\label{alg_plain_dw}
\algloop{select}
\algloop{repeat}
\begin{algorithmic}
\State Initialise the lower bound $l = -\infty$ and the upper bound $u = \infty$.
\State Solve the \gls{LPR} and use the solution as the initial values of $y$.
\For{$k$ in $\{1, 2, \ldots\}$}
	\State Solve the pricing subproblems.
	\State Compute the lower bound $\dualfunction(y)$ and update $l$.
  \State Update and solve the regularised \gls{RMP} and set $y$ to the solution.
	\State Run primal heuristics and update $u$.
\EndFor
\end{algorithmic}
\end{algorithm}

%% file: ph_sec_primal_heuristics.tex
In this section, we first review primal heuristics used with Dantzig-Wolfe decomposition or Lagrangian relaxation in the literature.
Then we propose two new primal heuristics: one based on decomposition and one based on machine learning.

\subsection{Review of primal heuristics based on decomposition}
\label{ph_subsec_review_of_primal_heuristics_based_on_decomposition}

\citet{MerlinAndSnadrin1983} solved \gls{UC} by applying Lagrangian relaxation and using a subgradient method for the dual problem.
In each iteration of the subgradient method, the pricing subproblems \eqref{eq_pricing_problem} were solved.
They then constructed a primal solution using the pricing subproblem solutions and tested if this was feasible for \eqref{problem_eq_mip}.
If the generation capacity was not sufficient to meet the demand or the spinning reserve at some time periods, they modified the dual step direction to increase the dual variable corresponding to the violated constraint to encourage more generators to be on at the shortage periods.
We note that this modification affects the optimisation of the dual variable and the resulting lower bounds may be suboptimal.

Instead of modifying the dual step direction, \citet{Guanetal1992} fixed infeasibility of the primal solutions by applying local search.
In each iteration, a primal solution was obtained from the subproblem solutions.
If this was infeasible, the cheapest available generators were committed to meet the demand and the spinning reserve.
After a feasible commitment decision was found, the amount of power output of each generator was optimised.
That is, the values of the binary variables in the original problem were fixed to the values corresponding to the feasible commitment and the resulting \gls{LP} was solved to compute the amount of power output.
In this paper, we refer to this approach as the {\em feasibility recovery local search primal heuristic}.

Another heuristic based on decomposition is the {\em column combination primal heuristic}.
In this heuristic, the solutions to the pricing subproblems are stored in a pool.
Then a constraint is added to the problem \eqref{problem_eq_mip} to restrict the pattern of the solution to those in the pool.
Let $\{ x_{\idxcomponent i}' = (u_{\idxcomponent i}, z_{\idxcomponent i}) \in \mathbb{R}^n \times \{0, 1\}^m \mid i \in J_{\idxcomponent} \}$ be the pool of solutions to pricing subproblem $\idxcomponent$ where $J_\idxcomponent$ is the index set.
Then, binary variables $w_{\idxcomponent i}$ ($\oneofidxcomponents$, $i \in J_\idxcomponent$) are added to the original problem \eqref{problem_eq_mip} together with the following constraints
\begin{gather*}
z_\idxcomponent = \sum_{i \in J_\idxcomponent} w_{\idxcomponent i} z_{\idxcomponent i}', \qquad \oneofidxcomponents, \\
\sum_{i \in J_\idxcomponent} w_{\idxcomponent i} = 1, \qquad \oneofidxcomponents, \\
w_{\idxcomponent i} \ge 0, \qquad \oneofidxcomponents, i \in J_\idxcomponent.
\end{gather*}
This is referred to as the restricted master IP by \citet{Vanderbeck2005}.
A standard \gls{MILP} solver can be used to solve the problem but the solution space is much smaller than the original problem.
This method was used to solve \gls{UC} by \citet{Takritietal2000} and to solve a stochastic version of \gls{UC} problems by \citet{Schulzeetal2017}.

\subsection{RMP partial-fixing}
\label{ph_subsec_primal_heuristics_based_on_rmp}

In this section, we introduce a new primal heuristic which uses the \gls{RMP} to construct primal feasible solutions.
In the following, we assume that the \gls{RMP} \eqref{problem_eq_rmp} is feasible, which may be ensured by adding some artificial columns.
First we solve the \gls{RMP} \eqref{problem_eq_rmp} without regularisation and for each subproblem $s$ compute a weighted-average of columns
\begin{equation*}
\hat{x}_\idxcomponent := \sum_{i \in \hat{I}_\idxcomponent} x_{\idxcomponent i} p_{\idxcomponent i}, \qquad \oneofidxcomponents,
\end{equation*}
where $\hat{I}_\idxcomponent$ is the index set of columns $\{x_{\idxcomponent i} \mid i \in \hat{I}_\idxcomponent\}$ used to formulate the \gls{RMP} for each $\idxcomponent$.
Using solutions $\hat{x}_\idxcomponent$ for each subproblem $\idxcomponent$, we construct a primal solution $\hat{x}$.
Although the integer variables in $x_{\idxcomponent i}$ satisfy the integrality constraint for any $\idxcomponent$ and $i \in \hat{I}_\idxcomponent$, those in $\hat{x}$ may not.
In this primal heuristic, we check whether the elements of $\hat{x}$ that correspond to binary decisions satisfies the integrality constraint and fix those which do.
In this way, we obtain a partially-fixed \gls{UC} problem, which is then solved by an \gls{MILP} solver.
In each iteration of the column generation procedure, we repeat the above process.
We refer this primal heuristic as {\em \gls{RMP} partial-fixing primal heuristic}.

Remarkably, the number of elements of $\hat{x}$ which violate the integrality constraint is bounded by a constant independent of the number of generators assuming the Simplex Method is used to solve the \gls{RMP}.
Thus, the difficulty to solve the partially-fixed \gls{MILP} is bounded.
This is a notable feature of this approach compared with typical primal heuristics.
For example, the problem solved in the column combination primal heuristic gets harder as the number of generator increases.

The above property can be shown using an argument similar to the one given by \citet{Bertsekasetal1983}.
Let $C$ be the number of complicating constraints in \eqref{problem_eq_mip}.
That is, $a \in \mathbb{R}^C$.
If we use the formulation shown in Appendix \ref{sec_appendix_problem_formulation}, $S$ is the number of generators and $C$ is double the number of time periods.
In the following we assume that $S > C$.
We note that in the \gls{RMP} \eqref{problem_eq_rmp}, there are $S + C$ equality constraints.
Thus, any basic solution to the \gls{RMP} has at most $S + C$ variables positive.
The second constraint in \eqref{problem_eq_rmp} ensures that for each $\idxcomponent$ at least one variable among $\{p_{\idxcomponent i}\}_{i \in \hat{I}_\idxcomponent}$ is positive.
Therefore, at most $C$ generators have two or more positive values among $\{p_{\idxcomponent i}\}_{i \in \hat{I}_\idxcomponent}$.
In other words, at least $S - C$ generators have exactly one positive value of $\{p_{\idxcomponent i}\}_{i \in \hat{I}_\idxcomponent}$.
In such a case, exactly one of $\{p_{\idxcomponent i}\}_{i \in \hat{I}_\idxcomponent}$ is equal to 1 and all of the others are equal to 0.
Thus, $\hat{x}_\idxcomponent$ equals exactly one of $\{x_{\idxcomponent i}\}_{i \in \hat{I}_\idxcomponent}$ and has integer values.

In practice the number of generators with multiple non-zeros in $\{p_{\idxcomponent i}\}_{i \in \hat{I}_\idxcomponent}$ may be smaller than $C$.
Furthermore, even if multiple elements in $\{p_{\idxcomponent i}\}_{i \in \hat{I}_\idxcomponent}$ are positive, only a small part of $\hat{x}_\idxcomponent$ may have fractional values.

We note that the \gls{RMP} \eqref{problem_eq_rmp} is not necessarily feasible.
Typically we need to run a few column generation iterations to gather enough columns to make the \gls{RMP} feasible.
Furthermore, even if the \gls{RMP} is feasible, the partially-fixed \gls{MILP} is not necessarily feasible.
In our experiments on \gls{UC} instances with practical data, we do not observe instances where the partially-fixed schedule is infeasible and typically the above primal heuristic successfully finds a primal feasible solution.
However, on some instances, it fails to provide a solution with small suboptimality, such as 0.1\%.
To handle such cases, we modify the method as follows.
Every time when we run the method, we record the upper bound.
If the \gls{RMP} is feasible but the upper bound is not improved for a prescribed number of the successive iterations (3 iterations in our implementation), we relax integer variables which are adjacent in time periods to those with fractional values.
For example, if the on-off status of generator $g$ on time $t$ is set free, we unfix the on-off status of the same generator of time $t-1$ and $t+1$.
By relaxing more binary variables, the partially-fixed \gls{UC} gets harder to solve but more likely to provide a better solution.

\subsection{Partial-fixing based on machine learning}
\label{ph_subsec_primal_heuristics_with_pretrained_model}

All of the primal heuristics discussed above use information gathered through the execution of the column generation procedure (or the subgradient methods).
In this section, we consider a primal heuristic based on machine learning.
We assume that the problem \eqref{problem_eq_mip} is to be solved repeatedly with different demand data $\omega$.

In the training phase, we sample $\omega$ (e.g.\ from historical data) and solve as many training instances as possible.
In this way, we obtain the data set for each sample $\omega$ and the corresponding optimal values of the binary variables.
Then, we use the data set to create a prediction model which takes $\omega$ as input and predicts the value of the binary variables.
We consider two alternatives: a neural network model and a nearest neighbour model.

The neural network model we consider in this paper is a feed-forward neural network \cite{Bishop2006}.
The model takes $\omega$ as input and outputs a vector each of whose elements is a predicted probability of the corresponding binary variables to be 1.
The neural network model is trained as a standard binary classification problem.

The nearest neighbour model is also considered as an alternative model to predict the values of binary variables.
When solving a new instance the nearest neighbour model compares the problem parameter $\omega$ with those in the training data set.
A prescribed number of the closest neighbours are selected and the average of the values of the binary variables are computed and used as the prediction of the probability.

The output of a prediction model can be used to find a feasible solution.
The simplest approach is to round the prediction to the nearest integer.
However, such a solution is usually infeasible when the problem is highly constrained.
Instead, as described below, we use the prediction to fix only a subset of binary variables so that the problem size is reduced.
Similar ideas have been explored by \citet{Xavieretal2020} and \citet{Wang2021}.

Pick a threshold value $\alpha \in (0.5, 1]$.
If the prediction is larger than $\alpha$ (smaller than $1 - \alpha$), fix the corresponding binary variables to 1 (0), and leave all the other variables unfixed.
Then the resulting partially-fixed \gls{MILP} is solved with an \gls{MILP} solver.

Choosing a suitable threshold value is a subtle task.
If we fix many binary variables the problem becomes small and can be solved quickly.
However, fixing too many variables may result in infeasibility or unacceptably large suboptimality.
On the other hand, fixing fewer variables results in a harder problem which takes longer to solve.

Instead of fixing the threshold value to a single value a priori, we try various values adaptively.
Namely, we first try a small threshold value and solve the partially-fixed \gls{MILP}.
If the resulting problem is infeasible, or if the resulting problem is solved, we try a larger threshold value.
We refer to the method based on neural network and nearest neighbour as {\em neural network partial-fixing primal heuristic} and {\em nearest neighbour partial-fixing primal heuristic} respectively.

In the context of Dantzig-Wolfe decomposition, the above method can be combined with the column generation procedure.
At the end of each iteration of the column generation procedure, we run one of the above primal heuristics and solve partially-fixed \gls{MILP}s using an optimisation solver.
We may impose a limit on the amount of time spent by the primal heuristic.
Namely, after the solver spends a certain amount of time, it is halted and the next iteration of the column generation procedure is executed.
In the next run of the primal heuristic, the solver is resumed from where it was interrupted in the previous iteration.

%% file: ph_sec_numerical_experiments.tex
In this section, the proposed methods are evaluated on large-scale \gls{UC} instances.
To assess the scalability, we consider 3 different problem sizes; problems with 200, 600 and 1,000 generators.
In all cases, the length of the planning horizon is 48 hours with a time resolution of 1 hour.
The generator data is based on \citet{Borghettietal2002}.
Since their sets of generators contain 200 generators at most, we combine multiple sets to create larger ones.
For example, to create a \gls{UC} instance with 1,000 generators, we combine 5 distinct 200-generator sets.
Each generator is unique and combination does not introduce symmetry.
The demand data is based on the historical demand data in the UK published by National Grid ESO.\footnote{https://www.nationalgrideso.com/}
A detailed description of the problem formulation is given in Appendix \ref{sec_appendix_problem_formulation}.

All methods are implemented in Python.
IBM ILOG CPLEX\footnote{https://www.ibm.com/products/ilog-cplex-optimization-studio} is used as the optimisation solver and PyTorch to implement the neural networks.  
The experiments are done on Intel\textsuperscript{\tiny\textregistered} Xeon\textsuperscript{\tiny\textregistered} E5-2670.

\subsection{Evaluation: Time to close gap with Dantzig-Wolfe decomposition}

\subsubsection{Experimental Setups}

In this experiment we used the primal heuristics alongside Dantzig-Wolfe decomposition and measured computational time to find a primal feasible solution and prove that its suboptimality was smaller than a prescribed tolerance.

To prepare the dataset used by the neural network and nearest neighbour partial-fixing, as many training instances as possible were solved to 0.25\% optimality in the training phase.
The training budget was 24 hours on 8 CPU cores.
To solve each training instance we used Dantzig-Wolfe decomposition (Algorithm \ref{alg_plain_dw}) with the feasibility recovery local search primal heuristic.
The number of solved instances are reported in Table \ref{tab_nn_training_summary}.

After the dataset was constructed, a neural network model was trained to predict the values of the binary variables.
We used a feed-forward neural network with 2 hidden-layers of 400 units per layer with ReLU activation function.
For simplicity, the time to train a neural network model, which is shown in Table \ref{tab_nn_training_summary}, is not included in the training budget of 24 hours.
When solving a test instance, we used the threshold values of 0.8, 0.9, 0.95, 0.99 and 1.0.

\begin{table}
\begin{center}  
\small
\caption{Statistics on the training of neural network models}
\label{tab_nn_training_summary}
\begin{tabular*}{9cm}{
@{\extracolsep{\fill}}lrrr}
\toprule
& 200 & 600 & 1000 \\
\midrule
number of training instance & 15,419 & 7,241 & 4,886 \\
time to train a model (s) & 757 & 1,514 & 1,837 \\
\bottomrule
\end{tabular*}
\end{center}  
\end{table}

The nearest neighbour method used the same dataset as the neural network model.
When solving a test instance, the parameter of the instance was compared with those of training instances, and the 50 closest instances (in Euclidean distance) were chosen to compute the average values of the binary variables.

To solve a test instances, Dantzig-Wolfe decomposition (Algorithm \ref{alg_plain_dw}) was used and the primal heuristics were run in the end of each column generation iteration.
The time limit of the primal heuristics was set to
\begin{equation}
\text{(primal heuristic time limit)} = \text{(time spent to solve pricing subproblems)} \cdot \left( \frac{k}{10} + 2\right),
\label{ph_eq_time_limit}
\end{equation}
where $k$ is the iteration number.
When the iteration number is large, the lower bound is typically tight and the upper bound provided by the primal heuristics is loose unless more time is allocated to it.
Thus, we allowed the primal heuristics to use more time in later iterations.
We note that some primal heuristics may stop before reaching the time limit.
For example, the partially-fixed \gls{MILP} solved in the \gls{RMP} partial-fixing primal heuristic was often solved before the timelimit was reached.

For comparison, we also ran CPLEX on the original \gls{MILP} problem without decomposition.

\subsubsection{Results}

In this experiment 40 test instances are solved.
The demand data to construct the test instances are sampled from a different year than those of the instances used to train the neural network and the nearest neighbour model.
Table \ref{ph_tab_dw_summary} shows the number of instances solved within the time limit of 20 minutes, the average computational time and the average number of column generation iterations to solve the instances or to reach the time limit.
When calculating the averages, the instances that are not solved within the time limit have the time of 20 minutes and the number of iterations reached by the 20-minute time limit is used.

\begin{table}
\begin{center}  
\caption{Computational time and required number of iterations}
\label{ph_tab_dw_summary}
\begin{tabular*}{\textwidth}{
@{\extracolsep{\fill}}l@{\hskip0.1cm}l@{\hskip0.1cm}
r@{\hskip0.1cm}r@{\hskip0.1cm}r
r@{\hskip0.1cm}r@{\hskip0.1cm}r
r@{\hskip0.1cm}r@{\hskip0.1cm}r}
\toprule
 & & \multicolumn{3}{l}{tol: 0.5\%} & \multicolumn{3}{l}{0.25\%} & \multicolumn{3}{l}{0.1\%} \\
\cmidrule{3-5}
\cmidrule{6-8}
\cmidrule{9-11}
size & method  & solved &  time & iter & solved &  time & iter & solved &  time & iter \\
\midrule
200  & feasibility recovery &       40 &           15.9 &         3.0 &       40 &           56.0 &        11.3 &        4 &           551.8 &        68.6 \\
     & column combination &       40 &           32.8 &         7.1 &       40 &           33.6 &         7.2 &       30 &           276.3 &        35.6 \\
     & RMP partial-fixing &       40 &           22.8 &         5.5 &       40 &           24.0 &         5.8 &       40 &   \textbf{37.5} &         8.4 \\
     & network partial-fixing &       40 &   \textbf{9.3} &         1.0 &       40 &  \textbf{12.0} &         1.3 &       40 &            56.3 &         5.2 \\
     & nearest partial-fixing &       40 &           11.6 &         1.2 &       40 &           15.5 &         1.6 &       40 &            74.4 &         6.3 \\
     & CPLEX &       40 &          215.8 &         0.0 &       37 &          336.0 &         0.0 &       18 &           558.0 &         0.0 \\
\midrule
600  & feasibility recovery &       40 &           42.2 &         2.0 &       40 &           84.4 &         5.4 &       30 &           263.5 &        16.8 \\
     & column combination &       40 &           83.9 &         5.2 &       40 &           86.1 &         5.3 &       40 &           105.6 &         6.4 \\
     & RMP partial-fixing &       40 &           60.7 &         4.4 &       40 &           62.8 &         4.6 &       40 &   \textbf{76.8} &         5.7 \\
     & network partial-fixing &       40 &  \textbf{39.6} &         1.3 &       40 &  \textbf{41.7} &         1.3 &       40 &           128.0 &         4.2 \\
     & nearest partial-fixing &       40 &           51.8 &         1.6 &       40 &           60.2 &         1.8 &       40 &           136.0 &         4.3 \\
     & CPLEX &        5 &          589.7 &         0.0 &        5 &          589.7 &         0.0 &        1 &           593.5 &         0.0 \\
\midrule
1000 & feasibility recovery &       40 &  \textbf{66.6} &         1.8 &       40 &          117.9 &         4.2 &       37 &           236.1 &         9.4 \\
     & column combination &       40 &          114.1 &         3.9 &       40 &          120.1 &         4.1 &       40 &           166.4 &         5.6 \\
     & RMP partial-fixing &       40 &           85.5 &         3.8 &       40 &           91.4 &         4.1 &       40 &  \textbf{113.7} &         5.2 \\
     & network partial-fixing &       40 &           70.8 &         1.3 &       40 &  \textbf{75.1} &         1.4 &       40 &           181.6 &         3.9 \\
     & nearest partial-fixing &       40 &           91.6 &         1.8 &       40 &           94.3 &         1.9 &       40 &           208.8 &         4.3 \\
     & CPLEX &        1 &          597.4 &         0.0 &        1 &          597.4 &         0.0 &        0 &           600.0 &         0.0 \\
\bottomrule
\end{tabular*}
\end{center}  
\end{table}

When the tolerance is loose, such as 0.5\% or 0.25\%, the neural network partial-fixing usually performs best.
We observe that the number of column generation iterations in these cases is very small and with averages all less than 2.
That is, on typical instances the neural network partial-fixing very quickly find a primal feasible solution satisfying the suboptimality tolerance with 0.25\% and Dantzig-Wolfe decomposition give a lower bound to assert that the suboptimality was smaller than 0.25\%.
The nearest neighbour partial-fixing is second best in half of the cases.
However, for all the cases, the average performance of the neural network partial-fixing is better than that of the nearest neighbour partial-fixing, in terms of the computational time and the number of column generation iterations.
The other primal heuristics are based on decomposition and require more column generation iterations to find primal feasible solutions of acceptable suboptimality.
This results in longer computational time for many instances.

If the target tolerance is tight (0.1\%), the results are different.
The \gls{RMP} partial-fixing outperforms the other methods in all the cases.
The neural network partial-fixing requires a smaller number of column generation iterations but needs longer computational time.
This is because the \gls{RMP} partial-fixing does not use all the time allocated to the primal heuristics but the neural network partial-fixing always run as long as the time limit.
The neural network and nearest neighbour partial-fixing primal heuristics found a primal solution with suboptimality smaller than 0.1\% for all the test instances.
However, the neural network partial-fixing is on average slower than the \gls{RMP} partial-fixing in the all cases and the nearest neighbour partial-fixing primal heuristic is even slower.
The column combination primal heuristic failed to find a primal solution for some test instances for the 200-generator case.
However, it successfully found a primal solution on all the test instances of size 600 and 1,000 and the average computational time is faster than the neural network partial-fixing primal heuristic.
The feasibility recovery local search primal heuristic also has better performance for larger test instances but is still inferior to the other primal heuristics.

\subsubsection{Analysis of the effect of training budget}

In this section, we study the effect of the training budget.
We consider cases where the training budget is 6, 12, 36 or 48 hours instead of 24 hours and observe how the performance of the methods is affected.
To this end, the neural network model is trained using the dataset generated within each of these training budgets.
The performance of the models is then evaluated as before and the results are reported in Table \ref{tab_training_budget}.
In all cases, all of the test instances are solved to within 0.1\% tolerance.

In many cases, both the neural network model and the nearest neighbour model tend to perform better with a larger training dataset.
Comparing the neural network models with 6-hour training, those with 48-hour training are all on average faster.
However, there is not a systematic improvement in the performance beyond 24-hour of training.
The room for additional performance gain seems limited if further training budget is given.
Similar discussion holds for the result of the nearest neighbour model.

\begin{table}
\caption{Performance of models with different training budgets}
\label{tab_training_budget}
\begin{tabular*}{\textwidth}{
@{\extracolsep{\fill}}
l@{\hskip0.1cm}l@{\hskip0.1cm}l@{\hskip0.1cm}
r@{\hskip0.1cm}r@{\hskip0.1cm}r
r@{\hskip0.1cm}r@{\hskip0.1cm}r
r@{\hskip0.1cm}r@{\hskip0.1cm}r}
\toprule
& & & \multicolumn{3}{l}{tol: 0.5\%} & \multicolumn{3}{l}{0.25\%} & \multicolumn{3}{l}{0.1\%} \\
\cmidrule{4-6}
\cmidrule{7-9}
\cmidrule{10-12}
size & method & budget &  solved &  time & iter &  solved &  time & iter &  solved &   time & iter \\
\midrule
200  & network & 6  &      40 &  7.7 &  1.1 &      40 &   9.4 &  1.7 &      40 &  30.5 &  6.8 \\
     &         & 12 &      40 &  7.8 &  1.1 &      40 &   8.9 &  1.5 &      40 &  29.5 &  6.6 \\
     &         & 24 &      40 &  7.5 &  1.1 &      40 &   9.3 &  1.6 &      40 &  27.8 &  6.3 \\
     &         & 36 &      40 &  7.5 &  1.1 &      40 &   8.6 &  1.4 &      40 &  25.4 &  6.0 \\
     &         & 48 &      40 &  7.3 &  1.0 &      40 &   8.4 &  1.4 &      40 &  26.3 &  6.3 \\
\midrule
     & neighbour & 6  &      40 &  8.0 &  1.2 &      40 &  10.2 &  1.9 &      40 &  41.9 &  8.6 \\
     &         & 12 &      40 &  7.7 &  1.1 &      40 &  10.7 &  2.1 &      40 &  40.8 &  8.6 \\
     &         & 24 &      40 &  7.6 &  1.1 &      40 &   9.8 &  1.8 &      40 &  38.2 &  8.2 \\
     &         & 36 &      40 &  7.6 &  1.1 &      40 &   9.3 &  1.6 &      40 &  34.0 &  7.6 \\
     &         & 48 &      40 &  7.9 &  1.2 &      40 &  10.2 &  1.9 &      40 &  31.4 &  7.3 \\
\midrule
600  & network & 6  &      40 &  31.2 &  1.5 &      40 &  33.7 &  1.8 &      40 &  71.2 &  5.3 \\
     &         & 12 &      40 &  27.8 &  1.2 &      40 &  32.1 &  1.5 &      40 &  65.8 &  4.8 \\
     &         & 24 &      40 &  27.5 &  1.1 &      40 &  30.3 &  1.4 &      40 &  64.6 &  4.7 \\
     &         & 36 &      40 &  29.7 &  1.4 &      40 &  30.6 &  1.4 &      40 &  64.0 &  4.7 \\
     &         & 48 &      40 &  30.2 &  1.4 &      40 &  30.2 &  1.4 &      40 &  64.5 &  4.7 \\
\midrule
     & neighbour & 6  &      40 &  33.6 &  1.8 &      40 &  37.4 &  2.1 &      40 &  83.7 &  6.0 \\
     &         & 12 &      40 &  32.7 &  1.6 &      40 &  36.8 &  2.0 &      40 &  75.4 &  5.6 \\
     &         & 24 &      40 &  32.0 &  1.5 &      40 &  35.0 &  1.9 &      40 &  76.1 &  5.5 \\
     &         & 36 &      40 &  32.4 &  1.6 &      40 &  35.3 &  1.9 &      40 &  71.7 &  5.4 \\
     &         & 48 &      40 &  34.2 &  1.8 &      40 &  37.3 &  2.1 &      40 &  73.8 &  5.6 \\
\midrule
1000 & network & 6  &      40 &  50.3 &  1.5 &      40 &  53.9 &  1.7 &      40 &  105.7 &  5.0 \\
     &         & 12 &      40 &  50.5 &  1.4 &      40 &  53.2 &  1.6 &      40 &   97.9 &  4.5 \\
     &         & 24 &      40 &  44.6 &  1.2 &      40 &  45.1 &  1.2 &      40 &   89.7 &  4.2 \\
     &         & 36 &      40 &  45.5 &  1.2 &      40 &  45.9 &  1.3 &      40 &   91.3 &  4.3 \\
     &         & 48 &      40 &  47.2 &  1.4 &      40 &  48.4 &  1.4 &      40 &   89.5 &  4.2 \\
\midrule
     & neighbour & 6  &      40 &  56.8 &  1.9 &      40 &  63.3 &  2.4 &      40 &  143.5 &  6.6 \\
     &         & 12 &      40 &  55.2 &  1.9 &      40 &  59.7 &  2.1 &      40 &  122.4 &  5.7 \\
     &         & 24 &      40 &  62.5 &  2.1 &      40 &  69.5 &  2.5 &      40 &  132.6 &  6.0 \\
     &         & 36 &      40 &  47.8 &  1.4 &      40 &  52.4 &  1.6 &      40 &  110.7 &  5.1 \\
     &         & 48 &      40 &  52.1 &  1.6 &      40 &  57.3 &  1.9 &      40 &  112.9 &  5.2 \\
\bottomrule
\end{tabular*}
\end{table}

\subsubsection{Analysis on neural network model architecture}

In the following, the effect of the model architecture is studied.
In the previous experiments, small feed-forward neural network models with 2 hidden-layers of 400 units per layer were considered.
Here, we additionally train deeper neural network models and measure their performances.
A deeper model consists of 4 hidden layers of 1000 units per layer.
We use the same training dataset which is generated with the training budget of 24 hours.
The performance of the models are evaluated similarly and the results are reported in Table \ref{tab_deeper_model}.
The difference in performance is relatively small.
Although we observed that the performance of the neural network model is noticeably better than the nearest neighbour model, there is no systematic advantage of using the deeper, more expressive model.

\begin{table}
\begin{center}  
\caption{Performance of the original and deeper neural network models}
\label{tab_deeper_model}
\begin{tabular*}{\textwidth}{
@{\extracolsep{\fill}}l@{\hskip0.1cm}l@{\hskip0.1cm}
r@{\hskip0.1cm}r@{\hskip0.1cm}r
r@{\hskip0.1cm}r@{\hskip0.1cm}r
r@{\hskip0.1cm}r@{\hskip0.1cm}r}
\toprule
 & & \multicolumn{3}{l}{tol: 0.5\%} & \multicolumn{3}{l}{0.25\%} & \multicolumn{3}{l}{0.1\%} \\
\cmidrule{3-5}
\cmidrule{6-8}
\cmidrule{9-11}
size & model  & solved &  time & iter & solved &  time & iter & solved &  time & iter \\
\midrule
200 & original &       40 &     7.5 &    1.0 &       40 &     9.3 &    1.6 &       40 &    27.8 &    6.4 \\
     & deeper &       40 &     7.4 &    1.0 &       40 &     8.6 &    1.4 &       40 &    28.5 &    6.6 \\
\midrule
600 & original &       40 &    27.5 &    1.2 &       40 &    30.3 &    1.4 &       40 &    64.6 &    4.7 \\
     & deeper &       40 &    28.6 &    1.2 &       40 &    33.3 &    1.6 &       40 &    63.1 &    4.6 \\
\midrule
1000 & original &       40 &    44.6 &    1.2 &       40 &    45.1 &    1.2 &       40 &    89.7 &    4.2 \\
     & deeper &       40 &    43.9 &    1.1 &       40 &    44.3 &    1.2 &       40 &    97.5 &    4.6 \\
\bottomrule
\end{tabular*}
\end{center}  
\end{table}

\subsection{Evaluation: Best upper bound within time limit}

\subsubsection{Setups}

All of the discussions so far aimed to find a primal feasible solution with guaranteed suboptimality (e.g.\ 0.1\%).
In this section, we consider a case where we are only interested in obtaining as good primal feasible solution as possible within a prescribed time budget.
This is of interest when we do not necessarily have enough time to achieve a given proven tolerance.

To evaluate the performance of the primal heuristics in this setup, we compare them by the quality (suboptimality) of feasible solutions found within a prescribed time limit.
The neural network and nearest neighbour partial-fixing primal heuristics do not require Dantzig-Wolfe decomposition to be run and so are used as stand-alone methods.
That is, we formulate the partially-fixed \gls{MILP} instances and solve them sequentially (from those with small threshold values).
We note that even though we are not interested in computing a lower bound, the other primal heuristics still require Dantzig-Wolfe decomposition to be run.

We use the same neural network and nearest neighbour models as in the previous experiments.
The configuration of the other primal heuristics are the same as well.

\subsubsection{Results}

For evaluation, the same 40 test instances are used.

The results are shown in Table \ref{ph_tab_standalone_summary}.
The columns labelled as 'solved' are the number of test instances where the primal heuristics found a feasible solution within the time limits of 1, 2 and 5 minutes respectively.
Among the instances where the primal heuristics found a feasible solution, the gap between the best upper bounds found by the primal heuristics and the best lower bounds found by running CPLEX for 4 hours are computed and shown in the columns labelled as 'suboptimality'.

When the time limit is 1 minute, just finding a primal feasible solution may not be trivial.
On the instances of size 200, all of the primal heuristics can find feasible solutions.
However, on larger instances such as those of size 600 or 1000, only the neural network partial-fixing and the feasibility recovery local search primal heuristic found a feasible solution on all of the test instances.
We also note that the neural network partial-fixing found primal feasible solutions of smaller suboptimality compared with the feasibility recovery local search primal heuristic on average.
The column combination primal heuristic failed to find any primal feasible solutions on more than half of the test instances of size 1000 within 1 minute.

With longer time limit, the primal heuristics are more likely to find primal feasible solutions.
On 200-generator case, the neural network partial-fixing performs best on average on any time limits.
However, on 600-generator and 1000-generator case, given sufficiently long time limit, such as 5 minutes, the \gls{RMP} partial-fixing finds the primal feasible solution of the smallest suboptimality among the other primal heuristics on average and the neural network partial-fixing finds the second best solutions.
We note that on any setups, the nearest neighbour partial-fixing performs worse than the neural network partial-fixing.

\begin{table}
\begin{center}  
\caption{Quality of feasible solutions found within time limits}
\label{ph_tab_standalone_summary}
\begin{tabular*}{\textwidth}{@{\extracolsep{\fill}} llrrrrrr}
\toprule
     &                        & \multicolumn{2}{l}{1 minute} & \multicolumn{2}{l}{2 minutes} & \multicolumn{2}{l}{5 minutes} \\
\cmidrule{3-4}
\cmidrule{5-6}
\cmidrule{7-8}
size & method                 & solved & subopt. & solved & subopt. & solved & subopt. \\
\midrule
200  & feasibility recovery &     40 &           0.208 &     40 &           0.159 &     40 &           0.147 \\
     & column combination &     40 &           0.094 &     40 &           0.084 &     40 &           0.078 \\
     & RMP partial-fixing &     40 &           0.056 &     40 &           0.050 &     40 &           0.047 \\
     & network partial-fixing &     40 &  \textbf{0.039} &     40 &  \textbf{0.034} &     40 &  \textbf{0.028} \\
     & nearest partial-fixing &     40 &           0.054 &     40 &           0.043 &     40 &           0.032 \\
\midrule
600  & feasibility recovery &     40 &           0.317 &     40 &           0.172 &     40 &           0.088 \\
     & column combination &     10 &           0.107 &     32 &           0.062 &     40 &           0.024 \\
     & RMP partial-fixing &     22 &           0.117 &     40 &  \textbf{0.028} &     40 &  \textbf{0.021} \\
     & network partial-fixing &     40 &  \textbf{0.067} &     40 &           0.039 &     40 &           0.026 \\
     & nearest partial-fixing &     39 &           0.242 &     40 &           0.046 &     40 &           0.032 \\
\midrule
1000 & feasibility recovery &     40 &           0.301 &     40 &           0.277 &     40 &           0.072 \\
     & column combination &      0 &              - &     23 &           0.112 &     40 &           0.033 \\
     & RMP partial-fixing &      5 &           0.220 &     40 &           0.079 &     40 &  \textbf{0.014} \\
     & network partial-fixing &     40 &  \textbf{0.243} &     40 &  \textbf{0.042} &     40 &           0.028 \\
     & nearest partial-fixing &     39 &           0.550 &     40 &           0.052 &     40 &           0.035 \\
\bottomrule
\end{tabular*}
\end{center}  
\end{table}

%% file: ph_sec_conclusion.tex
We have considered primal heuristics based on decomposition and based on machine learning to solve the \gls{UC} problem together with Dantzig-Wolfe decomposition.
The primal heuristic based on decomposition uses the \gls{RMP} to fix a subset of the integer variables.
The \gls{MILP} solver could quickly solve the resulting partially-fixed \gls{MILP} and provide feasible solutions with small suboptimality, such as 0.1\%.
We discussed that the number of unfixed integer variables were bounded by a constant independent of the number of generators, which may explain to some extent the scalability of our method to solve large-scale instances.
On the other hand, the other primal heuristic uses a pre-trained neural network.
Given the prediction of a pre-trained model on the values of binary variables, we proposed an approach to reduce the problem size by fixing some of the binary variables.
By adaptively adjusting the number of variables to be fixed, the method can find a good solution quickly.
Furthermore, with Dantzig-Wolfe decomposition, our approach can find a solution of proven suboptimality.

In our numerical experiments, we compared our methods with primal heuristics from the literature, the feasibility recovery local search primal heuristic and the column combination primal heuristic.
Furthermore, we implemented another version of our primal heuristic which used a nearest neighbour approach instead of the neural network.
All of our methods were able to find a primal feasible solution of 0.1\% suboptimality on all of the test instances, while the primal heuristics from the literature sometimes failed.
The experiments also revealed that the primal heuristic based on machine learning was superior when the suboptimality tolerance was relatively large, such as 0.5\% or 0.25\%, while the primal heuristic based on decomposition was the best when the tolerance was small, for example 0.1\%.
In any cases the primal heuristic based on the neural network outperformed the nearest neighbour approach.

We also considered a case where we were only interested in the best primal solution found in a prescribed time budget.
In this case, the primal heuristic based on machine learning was used as a standalone method without Dantzig-Wolfe decomposition and outperformed other primal heuristics when the timelimit is tight (2 minutes or less).
If the timelimit is long, such as 5 minutes, the primal heuristic based on the \gls{RMP} often performed better.
Similar to the previous experiment, we observed that the primal heuristic based on the neural network was superior to the nearest neighbour approach in all cases.

%% file: ph_appendix.tex
\appendix

\section{Problem formulation}
\label{sec_appendix_problem_formulation}

We follow one of the standard formulations in literature and formulate the following constraints:
\begin{itemize}
\item \textbf{Load balance}: Generators have to meet all the demand in each time period (generation shedding of 0 cost is allowed).
\item \textbf{Reserve}: To deal with contingencies, it is required to keep a sufficient amount of back up in each time period, which can be activated quickly.
\item \textbf{Power output bounds}: Each generator's power output has to be within its limit.
\item \textbf{Ramp rate}: Generators can only change their outputs within the ramp rates.
\item \textbf{Minimum up/downtime}: If switched on (off), each generator has to stay on (off) for a given minimum period.
This is to avoid thermal stress in the generators which may cause wear and tear of the turbines.
\end{itemize}

The formulation of the model is as follows.

\begin{itemize}
\item{Parameters}
\begin{itemize}
\item $G$: The number of generators
\item $T$: The number of time periods where decisions are taken
\item $C^{\mathrm{nl}}_{g}$: no-load cost of generator $g$
\item $C^{\mathrm{mr}}_{g}$: marginal cost of generator $g$
\item $C^{\mathrm{up}}_{g}$: startup cost of generator $g$
\item $P^{\max/\min}_{g}$: maximum/minimum generation limit of generator $g$
\item $P^{\mathrm{ru}/\mathrm{rd}}_{g}$: operating ramp up/down limits of generator $g$
\item $P^{\mathrm{su}/\mathrm{sd}}_{g}$: startup/shutdown ramp limits of generator $g$
\item $T^{\mathrm{u}/\mathrm{d}}_{g}$: minimum uptime/downtime of generator $g$
\item $P^{\mathrm{d}}_{t}$: power demand at time $t$
\item $P^{\mathrm{r}}_{t}$: reserve requirement at time $t$
\end{itemize}

\item{Variables}
\begin{itemize}
\item $\alpha_{gt} \in \{0, 1\}$: 1 if generator $g$ is on in period $t$, and 0 otherwise
\item $\gamma_{gt} \in \{0, 1\}$: 1 if generator $g$ starts up in period $t$, and 0 otherwise
\item $\eta_{gt} \in \{0, 1\}$: 1 if generator $g$ shut down in period $t$, and 0 otherwise
\item $p_{gt} \ge 0$: power output of generator $g$ in period $t$
\end{itemize}

\item The objective is the total cost
\[
\min \sum_{t = 1}^T \sum_{g = 1}^G
\left( C^{\mathrm{nl}}_g \alpha_{gt} + C^{\mathrm{mr}}_g p_{gt} +
C^{\mathrm{up}}_g \gamma_{gt}
\right).
\]
This is to be minimised subject to the following constraints.

\item Load balance
\begin{equation*}
\sum_{g = 1}^G p_{gt} \ge P^{\mathrm{d}}_{t}
\qquad t = 1, 2, \ldots, T.
\label{eq:uc_first_constraint}
\end{equation*}

\item Reserve
\begin{equation*}
\sum_{g = 1}^G (P^{\max}_{g} \alpha_{gt} - p_{gt})
\ge P^{\mathrm{r}}_t
\qquad t = 1, 2, \ldots, T.
\end{equation*}

\item Power output bounds
\begin{equation*}
P^{\min}_{g} \alpha_{gt} \le p_{gt} \le P^{\max}_{g} \alpha_{gt}
\qquad g = 1, 2, \ldots, G, t = 1, 2, \ldots, T
\end{equation*}

\item Ramp rate
\begin{equation*}
p_{gt} - p_{g \, t-1} \le P^{\mathrm{ru}}_g \alpha_{g \, t-1}
+ P^{\mathrm{su}}_g \gamma_{gt}
\qquad g = 1, 2, \ldots, G, t = 2, 3, \ldots, T.
\end{equation*}
\begin{equation*}
p_{g \, t-1} - p_{gt} \le P^{\mathrm{rd}}_g \alpha_{gt}
+ P^{\mathrm{sd}}_g \eta_{gt}
\qquad g = 1, 2, \ldots, G, t = 2, 3, \ldots, T.
\end{equation*}

\item Minimum up/downtime
\begin{equation*}
\sum_{u=\max\{t-T^\mathrm{u}_g+1, 1\}}^t \gamma_{gu} \le \alpha_{gt}
\qquad g = 1, 2, \ldots, G, t = 1, 2, \ldots, T
\end{equation*}
\begin{equation*}
\sum_{u=\max\{t-T^\mathrm{u}_g+1, 1\}}^t \eta_{gu} \le 1 - \alpha_{gt}
\qquad g = 1, 2, \ldots, G, t = 1, 2, \ldots, T
\end{equation*}

\item Polyhedral/Switching constraints (to enforce binaries to work as we expect)
\begin{equation*}
\alpha_{gt} - \alpha_{g \, t-1} = \gamma_{gt} - \eta_{gt}
\qquad g = 1, 2, \ldots, G, t = 2, 3, \ldots, T
\end{equation*}
\begin{equation*}
1 \ge \gamma_{gt} + \eta_{gt}
\qquad g = 1, 2, \ldots, G, t = 2, 3, \ldots, T
\label{eq:uc_last_constraint}
\end{equation*}
\end{itemize}

We note that the complicating constraints are inequality in the above formulation but the discussion in this paper (e.g.\ the number of fractional values in the solution to the \gls{RMP}) still holds.